\documentclass[10pt]{article}
\usepackage[english]{babel}
\usepackage{amsfonts}
\usepackage[dvips]{graphicx}

\begin{document}

\title{\bf On local symbolic approximation and resolution of ODEs using Implicit Function Theorem}

\date{January 2006}

\author{\bf Gianluca Argentini \\
\normalsize gianluca.argentini@riellogroup.com \\
\textit{Research \& Development Department}\\
\textit{Riello Burners}, 37048 San Pietro di Legnago (Verona), Italy}

\maketitle

\begin{abstract}
In this work the implicit function theorem is used for searching local symbolic resolution of differential equations. General results of existence for first order equations are proven and some examples, one relative to cavitation in a fluid, are developed. These examples seem to show that local approximation of non linear differential equations can give useful informations about symbolic form of possible solutions, and in the case a global solution is known, locally the accuracy of approximation can be good. 
\end{abstract}

\section{Introduction}

The use of {\it Implicit Function} or {\it Dini}'s Theorem for local resolution of differential equations is not new (see \cite{mingari} for an historical example of application due to Bernoulli, and \cite{valent} for an application of theoretical existence for a non linear system of PDEs). But in the technical literature this method is very little widespread and limited to some specific situations. Moreover, for some physical and engineering aims, often usual numerical resolution of differential equations is not very useful, because the knowledge of an analytical dependence from physical or geometrical parameters is required (see e.g \cite{buttazzo}) for the problem of searching analytical dependence for resistence coefficients in flows).\\
In this work some statements about local symbolic resolution of ordinary differential equations are established, and some examples, two very simple but illustrative the technique, and other two more technical and useful, are developed.

\section{ODEs and Dini's Theorem}

Let $[a,b]$, $a < b$, an interval of $\mathbb{R}$. If we denote with $y=y(x)$ a $C^n([a,b])$ function, for some $n$ integer, $n > 0$, defined on $[a,b]$ and with $y^{(k)}=y^{(k)}(x)$ its $k$-th derivative, an ordinary differential equation (ODE) is an expression like this:

\begin{equation}\label{ode1}
	F(x,y(x),y^{(1)}(x),...,y^{(n)(x)})=0 \hspace{0.5cm} \forall x \in [a,b]
\end{equation}

\noindent where $F$ is a function $F:[a,b] \times \mathbb{R}^{n+1} \rightarrow \mathbb{R}$. Let $\partial_x F$ the partial derivative of $F$ respect the $x$ variable in $[a,b]$ and $\partial_k F$, $0 \leq k \leq n$, the partial derivative respect a suitable variable of $\mathbb{R}^{n+1}$ domain.\\
Suppose now $F \in C^m\left([a,b] \times \mathbb{R}^{n+1}\right)$ for some integer $m > 0$.\\
If $T_0=(x_0^0,p_0^0,...,p_n^0) \in [a,b] \times \mathbb{R}^{n+1}$ is a point such that $F(T_0)=0$, and $\partial_k F(T_0) \neq 0$ for a particular value of $k$, then from $Dini$'s or {\it Implicit Function} Theorem there are a constant $\delta > 0$, a costant $M > 0$, a neighbour $U_0$ of $(p_0^0,...,p_n^0)$ and a unique function $\phi: [a,b] \times \mathbb{R}^{n} \rightarrow \mathbb{R}$, $\phi \in C^m\left([a,b]\right)$, such that

\begin{eqnarray}\label{dini1}
	F\left(x,p_0,...,p_{k-1},\phi(x,p,...,p_{k-1},p_{k+1},...,p_n),p_{k+1},...,p_n\right)=0 \\
	\nonumber{\forall (x,p_0,...,p_{k-1},p_{k+1},...,p_n) \in (x_0-\delta,x_0+\delta) \times U_0}
\end{eqnarray}

\begin{equation}\label{dini2}
	\phi(x_0^0,p_0^0,...,p_{k-1}^0,p_{k+1}^0,...,p_n^0)=p_k^0 \hspace{0.5cm}
\end{equation}

\begin{equation}\label{dini3}
	\left|\phi(x_0,p_0,...,p_{k-1},p_{k+1},...,p_n)-p_k^0\right| \leq M
\end{equation}

\noindent and

\begin{equation}\label{dini4}
	\frac{\partial \phi}{\partial p_i}=-\frac{\partial_i F}{\partial_k F}
\end{equation}

\noindent For a complete proof, see e.g. \cite{courant}.\\
Let we suppose $n = 1$, so the differential equation is $F(x,y(x),y'(x))=0$. If the equation has order greater then $1$, it is possible to write it as a system of order $1$ differential equations (see e.g. \cite{sanchez}), and Dini Theorem has a natural extension to systems of equations.

\newtheorem{theo}{Theorem}
\newtheorem{prop}{Proposition}

\begin{theo}\label{theo1}
	Let $F(x,y(x),y'(x))=0$ an {\normalfont ODE}, with $F \in C^1\left([a,b] \times \mathbb{R}^2\right)$ a real function $F(x,p,q)$, and let $(x_0,p_0,q_0)$ such that $x_0 \in (a,b)$, $F(x_0,p_0,q_0)=0$ and $\partial_q F(x_0,p_0,q_0) \neq 0$. Then exist $\delta > 0$, $\rho > 0$, a function $u \in C^1([x_0-\delta,x_0+\delta])$ and a function $\phi \in C^0([x_0-\delta,x_0+\delta] \times [p_0-\rho,p_0+\rho])$ such that $\phi(x_0,p_0)=q_0$ and 
	\begin{equation}\label{theoFormula1}
			u'(x)=\phi(x,u(x)), \hspace{0.2cm} F(x,u(x),\phi(x,u(x)))=0 \hspace{0.5cm} \forall x \in [x_0-\delta,x_0+\delta]
	\end{equation}
\end{theo}

\noindent {\it Proof}. From Dini Theorem exist $\delta_1 > 0$, $\rho_1 > 0$ and a function $\phi \in C^0([x_0-\delta_1,x_0+\delta_1] \times [p_0-\rho_1,p_0+\rho_1])$ such that 
\begin{equation}\label{F0temp}
	F(x,p,\phi(x,p))=0 \hspace{0.3cm} \forall (x,p) \in [x_0-\delta_1,x_0+\delta_1] \times [p_0-\rho_1,p_0+\rho_1]
\end{equation} 
\noindent We can consider the {\it differential problem}
\begin{equation}\label{ivproblem1}
	y'(x)=\phi(x,y(x)), \hspace{0.5cm} y(x_0)=p_0
\end{equation}
\noindent Being $\phi$ a continuous function, from the theorem of {\it Peano} about existence of solutions for first order ordinary differential problems (or by local existence theorems) (see \cite{brota},\cite{sanchez}), exist $\delta_2$ with $0 < \delta_2 \leq \delta_1$ and a function $u \in C^1([x_0-\delta_2,x_0+\delta_2])$ (local) solution of (\ref{ivproblem1}). Being $u$ a continuous function and $u(x_0)=p_0$, exists a $0 < \delta_3 \leq \delta_2$ such that $u(x) \in [p_0-\rho_1,p_0+\rho_1] \hspace{0.3cm} \forall x \in [x_0-\delta_3,x_0+\delta_3])$. Then, from (\ref{F0temp}), if $x \in [x_0-\delta_3,x_0+\delta_3]) \subseteq [x_0-\delta_1,x_0+\delta_1])$, we have $F(x,u(x),\phi(x,u(x)))=0$, so with the choice $\delta=\delta_3$ and $\rho=\rho_1$ the proof is completed. $\square$ \\

\noindent Therefore, under suitable hypothesis, ODE $F(x,y(x),y'(x))=0$ has at least local solutions. But the theorem says too that the original ODE can be {\it locally} written in a {\it new form}, i.e. in the expression $y'(x)=\phi(x,y(x))$, that is the {\it normal} form. In this sense, the term locally means that there is a neighbour of $(x_0,p_0,q_0)$ where the original ODE can be written in the new form, and a local solution $u$ exists such that $u(x_0)=p_0$ and $u'(x_0)=q_0$.\\

\noindent {\it Example 1}. Consider the ODE
	\begin{equation}
 		2y(x)-y'(x)=0
	\end{equation}
\noindent in the interval $[a,b]$ where $a<0$ and $b>0$. The function $F$ is $F(x,p,q)=2p-q$, with $F \in C^\infty \left([a,b] \times \mathbb{R}^2\right)$ and $\partial_q F = -1$. Let $x_0=0$ and $p_0=1$; then $q_0=2$. From theorem 1 exist two functions $u$ and $\phi$ such that $\phi(0,1)=2$ and $2p-\phi(x,p)=0$ in a neighbour of $(0,1)$. Then $\phi(x,p)=2p$, so from (\ref{theoFormula1}) $u'(x)=2u(x)$, that is $u(x)=ce^{2x}$. The constant $c$ is determined from the condition $u(0)=p_0=1$, hence $c=1$.\\

In the case $\partial_qF(x_0,p_0,q_0)=0$, we can consider the possibility of a local expression like $y(x)=\psi(x,y'(x))$, as stated in the following

\begin{theo}\label{theo2}
	Let $F(x,y(x),y'(x))=0$ an {\normalfont ODE}, with $F \in C^1\left([a,b] \times \mathbb{R}^2\right)$ a real function $F(x,p,q)$, and let $(x_0,p_0,q_0)$ such that $x_0 \in (a,b)$, $F(x_0,p_0,q_0)=0$ and $\partial_p F(x_0,p_0,q_0) \neq 0$. Then exist $\delta > 0$, $\rho > 0$ and a function $\psi \in C^0([x_0-\delta,x_0+\delta] \times [q_0-\rho,q_0+\rho])$ such that $\psi(x_0,q_0)=p_0$ and 
	\begin{equation}\label{theoFormula2}
			F(x,\psi(x,q),q)=0 \hspace{0.5cm} \forall (x,q) \in [x_0-\delta,x_0+\delta] \times [q_0-\rho,q_0+\rho]
	\end{equation}
If exists a solution $u \in C^1([x_0-\sigma,x_0+\sigma]$, with $0 < \sigma < \delta$, for the problem
	\begin{equation}\label{theo2IVproblem}
		\left\{
		\begin{array}{ll}
			$$y(x) = \psi(x,y'(x)) \hspace{0.5cm} x \in [x_0-\sigma,x_0+\sigma] $$\\
			$$y'(x_0) = q_0$$
		\end{array}
		\right.
	\end{equation}
\noindent then
	\begin{equation}\label{theo2sol}
		F(x,u(x),u'(x))=0 \hspace{0.5cm} \forall x \in [x_0-\sigma,x_0+\sigma]
	\end{equation}
\end{theo}

\noindent {\it Proof}. The first statement is a direct consequence of Dini's Theorem. Let then $u$ a solution of the problem (\ref{theo2IVproblem}), with $\sigma_1 < \delta$. From the fact that $u'$ is continuous in the interval $[x_0-\sigma_1,x_0+\sigma_1]$, we can choice the constant $\sigma$ so that $u'(x) \in [q_0-\rho,q_0+\rho]$ for every $x \in [x_0-\sigma,x_0+\sigma]$. Then, from (\ref{theoFormula2}) follows that $F(x,u(x),u'(x))=0$ for every $x \in [x_0-\sigma,x_0+\sigma]$. $\square$ \\

\noindent Therefore, a local solution for the new differential equation $y(x)=\psi(x,y'(x))$ is a local solution for the original ODE $F(x,y(x),y'(x))=0$.\\

\noindent {\it Example 2}. Consider the non linear ODE 
\begin{equation}
	y'(x)^2+y(x)^2-1=0
\end{equation}
\noindent The function $F$ is $F(x,p,q)=q^2+p^2-1$, with $\partial_pF=2p$ and $\partial_qF=2q$. If we consider the point $T_0=(\frac{\pi}{2},1,0)$, we have $\partial_qF(T_0)=0$ and $\partial_pF(T_0)=2$, so from the first statement of theorem 2, if $\left|q\right| \leq 1$ exists the function $\psi=\psi(x,q)=\sqrt{1-q^2}$ such that (locally) $p=\sqrt{1-q^2}$. Then we can consider the problem
	\begin{equation}
		\left\{
		\begin{array}{ll}
			$$y(x) = \sqrt{1-y'(x)^2} \hspace{0.5cm} \forall x \in [0,\pi] $$\\
			$$y'(\frac{\pi}{2}) = 0$$
		\end{array}
		\right.
	\end{equation}
\noindent The function $u(x)=sin(x)$ is a solution of this problem, with $\left|u'(x)\right| \leq 1$ as required, and is a local (but in this case global too) solution of the original ODE. Note that finding the solution of the original equation implies the use of complex variables: the internal command {\ttfamily DSolve[y'[x]\^{}2+y[x]\^{}2-1==0,y[x],x]} of the software {\it Mathematica} gives the integrals
\begin{equation}
	u(x)=1-2sin^2\left[\frac{1}{2}(x-ic)\right], \hspace{0.5cm} v(x)=1-2sin^2\left[\frac{1}{2}(-x-ic)\right]
\end{equation}
\noindent where $i$ is the imaginary unit and $c$ a complex constant. We obtain our previuos local solution respectively for, e.g., $c=i\frac{3}{2}\pi$ and $c=-i\frac{3}{2}\pi$.

\section{Local symbolic approximation for an ODE}

From implicit function theorem, if $F$ has continuous partial derivatives, using Taylor expansion near the point $(x_0,p_0)$ (or $(x_0,q_0)$) and the chain rule for the derivatives of composite function it is possibile to write the following expression for the function $\phi$ (or $\psi$):

\begin{eqnarray}\label{taylor1}
	\phi(x,p) = \phi(x_0,p_0) + \left[\partial_x \phi(x_0,p_0)\right](x-x_0) + \left[\partial_p \phi(x_0,p_0)\right](p-p_0) +\\
	\nonumber{+ \left[\frac{1}{2}\partial_{xx}^2 \phi(x_0,p_0)\right](x-x_0)^2 + \left[\frac{1}{2}\partial_{pp}^2 \phi(x_0,p_0)\right](p-p_0)^2 + ...}
\end{eqnarray}

\noindent with, e.g. for $x$-derivation,

\begin{eqnarray}\label{taylorx}
	\partial_x \phi(x_0,p_0) &=& -\frac{\partial_x F}{\partial_q F}(x_0,p_0,q_0) \\
	\partial_{xx}^2 \phi(x_0,p_0) &=& -\frac{\left(\partial_x \phi \right)^2 \partial_{qq}^2F + 	2\partial_x \phi \hspace{0.1cm}\partial_{xq}^2F + \partial_{xx}^2F}{\partial_q F} \label{taylorx2}
\end{eqnarray}

\noindent where the derivatives of $\phi$ are evaluated at $(x_0,p_0)$ and those of $F$ at $(x_0,p_0,q_0)$.\\
The expansion gives a local analytical approximation of the original ODE. If for this equation a global solution doesn't exist or it is not given by elementary functions, the expansion could be a more simple equation to solve. From previous theorems, a solution of this equation is an approximate local solution of the original ODE.\\

If $F \in C^1\left([a,b] \times \mathbb{R}^2\right)$, $F(x_0,p_0,q_0)=0$ and $\partial_q F(x_0,p_0,q_0) \neq 0$, from previous identities the first order differential equation $F(x,y(x),y'(x))$ has the following local {\it approximation of first order}:

\begin{eqnarray}\label{approxODE1}
	y'(x) = q_0 + \left[\frac{\partial_x F}{\partial_q F}(x_0,p_0,q_0)\right]x_0 + \left[\frac{\partial_p F}{\partial_q F}(x_0,p_0,q_0)\right]p_0 - \\ 
	\nonumber{- \left[\frac{\partial_x F}{\partial_q F}(x_0,p_0,q_0)\right]x - \left[\frac{\partial_p F}{\partial_q F}(x_0,p_0,q_0)\right]y(x)}
\end{eqnarray}

\noindent If $\partial_p F(x_0,p_0,q_0)=0$, the solution of this first order ODE is 

\begin{equation}
	y(x) = c + \left[q_0 + \frac{\partial_x F}{\partial_q F}x_0\right]x - \frac{1}{2}\frac{\partial_x F}{\partial_q F}x^2
\end{equation}

\noindent otherwise its general solution is (see e.g \cite{boyce})

\begin{equation}\label{approxSOL1}
	y(x) = p_0 + \frac{\partial_xF\partial_qF}{\left[\partial_pF\right]^2} + \frac{\partial_xFx_0 + \partial_qFq_0}{\partial_pF} - \frac{\partial_xF}{\partial_pF}x + c\hspace{0.1cm}e^{-\frac{\partial_pF}{\partial_qF}x}
\end{equation}

\noindent where $c$ is a constant and the partial derivatives are all calculated at point $(x_0,p_0,q_0)$. Imposing the condition $y(x_0)=p_0$ we obtain

\begin{equation}\label{constant1}
	c = - \frac{\left[\partial_xF + \partial_pFq_0\right]\partial_qF\hspace{0.1cm}e^{\frac{\partial_pF}{\partial_qF}x_0}}{\left[\partial_pF\right]^2}
\end{equation}

\noindent while, with this value of the constant $c$, we have $y'(x_0)=q_0$ as required.\\
So we have proven the following

\begin{prop}\label{prop1}
	In the hypothesis of Theorem 1, the {\normalfont ODE} $F(x,y(x),y'(x))=0$ has a local solution, approximated to first order, of the form
	\begin{eqnarray}
		y(x)&=&a_0 + b_0x + c_0x^2 \hspace{0.8cm} if \hspace{0.2cm} \partial_pF(x_0,p_0,q_0)=0 \\
		y(x)&=&a_0 + b_0x + c_0e^{d_0x} \hspace{0.5cm} if \hspace{0.2cm} \partial_pF(x_0,p_0,q_0) \neq 0
	\end{eqnarray}
\noindent where $a_0$, $b_0$, $c_0$ and $d_0$ are constants depending from $(x_0,p_0,q_0)$.
\end{prop}

\noindent {\it Example 1 bis}. For the ODE $2y(x)-y'(x)=0$, being $F(x,p,q)=2p-q$, at point $(0,1,2)$ we obtain $\partial_xF=0$, $\partial_pF=2$ and $\partial_qF=-1$, so that $a_0=b_0=0$, $c_0=1$ and $d_0=2$, therefore from Proposition 1 the local approximated solution of first order is $y(x)=e^{2x}$, which is the exact local (global) solution of the original equation, satisfying the condition $y(0)=1$.\\

If $F \in C^1\left([a,b] \times \mathbb{R}^2\right)$, $F(x_0,p_0,q_0)=0$ and $\partial_p F(x_0,p_0,q_0) \neq 0$, in a similar fashion the first order differential equation $F(x,y(x),y'(x))$ has the following local {\it approximation of first order}:

\begin{eqnarray}\label{approxODE2}
	y(x) = p_0 + \left[\frac{\partial_x F}{\partial_p F}(x_0,p_0,q_0)\right]x_0 + \left[\frac{\partial_q F}{\partial_p F}(x_0,p_0,q_0)\right]q_0 - \\ 
	\nonumber{- \left[\frac{\partial_x F}{\partial_p F}(x_0,p_0,q_0)\right]x - \left[\frac{\partial_q F}{\partial_p F}(x_0,p_0,q_0)\right]y'(x)}
\end{eqnarray}

\noindent If $\partial_qF(x_0,p_0,q_0)=0$, the local approximated solution is simply

\begin{equation}\label{approxSOL2simply}
	y(x) = p_0 + \frac{\partial_xF}{\partial_pF}x_0 - \frac{\partial_xF}{\partial_pF}x
\end{equation}

\noindent while in the case $\partial_qF(x_0,p_0,q_0) \neq 0$ one can verify that the local solution is the same as (\ref{approxSOL1}). So the following statement holds:

\begin{prop}\label{prop2}
	In the hypothesis of Theorem 2, the {\normalfont ODE} $F(x,y(x),y'(x))=0$ has a local solution, approximated to first order, of the form
	\begin{eqnarray}
		y(x) &=& a_0 + b_0x \hspace{1.9cm} if \hspace{0.2cm} \partial_qF(x_0,p_0,q_0)=0 \\
		y(x) &=& a_0 + b_0x + c_0e^{d_0x} \hspace{0.5cm} if \hspace{0.2cm} \partial_qF(x_0,p_0,q_0) \neq 0
	\end{eqnarray}
\noindent where $a_0$, $b_0$, $c_0$ and $d_0$ are constants depending from $(x_0,p_0,q_0)$.
\end{prop}
	
\noindent {\it Example 2 bis}. For the ODE $y'(x)^2+y(x)^2-1=0$, being $F(x,p,q)=p^2+q^2-1$, at point $(\frac{\pi}{2},1,0)$ we have $\partial_xF=0$, $\partial_pF=2$ and $\partial_qF=0$, therefore $a_0=1$, $b_0=0$, $c_0=0$ and $d_0=2$. From previous proposition, a local approximated solution of first order is $y(x)=1$. This result, compared with the exact solution $y(x)=sin(x)$, is not satisfactory. A greater order of approximation is required.\\

If $F \in C^n\left([a,b] \times \mathbb{R}^2\right)$ for a suitable integer $n > 0$, $F(x_0,p_0,q_0)=0$ and, e.g., $\partial_p F(x_0,p_0,q_0) \neq 0$, it is possible consider higher order derivatives for the implicit function $\psi$, as (\ref{taylorx2}) or

\begin{eqnarray}\label{taylorq2}
	\partial_{qq}^2 \psi(x_0,q_0) = -\frac{\left(\partial_q \psi \right)^2 \partial_{pp}^2F + 	2\partial_q \psi \hspace{0.1cm}\partial_{pq}^2F + \partial_{qq}^2F}{\partial_p F}
\end{eqnarray}

\noindent Then, from Taylor expansion of $\psi$, we can write the following approximation of order $(1,2)$ (1 for $x$-variable, 2 for $q$-variable) for an ODE:

\begin{eqnarray}\label{approxODE12}
	&&y(x) = p_0 + \nonumber \\
	&&+ \frac{1}{2\partial_pF}\left[ 2\partial_xFx_0 + 2\partial_qFq_0 - \left( \frac{\left[\partial_qF\right]^2\partial_{pp}^2F}{\left[\partial_pF\right]^2} - 2\frac{\partial_qF}{\partial_pF}\partial_{pq}^2F + \partial_{qq}^2F\right)q_0^2 \right] - \nonumber \\
	&&- \frac{1}{\partial_pF}\left[\partial_xFx + \left(\partial_qFq - \left( \frac{\left[\partial_qF\right]^2\partial_{pp}^2F}{\left[\partial_pF\right]^2} + 2\frac{\partial_qF}{\partial_pF}\partial_{pq}^2F + \partial_{qq}^2F\right)q_0\right)y'(x) \right] \nonumber \\
	&&- \frac{1}{2\partial_pF}\left( \frac{\left[\partial_qF\right]^2\partial_{pp}^2F}{\left[\partial_pF\right]^2} - 2\frac{\partial_qF}{\partial_pF}\partial_{pq}^2F + \partial_{qq}^2F\right)\left[y'(x)\right]^2 \nonumber
\end{eqnarray}

\noindent Expressions like this gives more accurate local solutions.\\

\noindent {\it Example 2 ter}. For the differential equation of Example 2 we have $\partial_xF(\frac{\pi}{2},1,0)=\partial_qF(\frac{\pi}{2},1,0)=\partial_{pq}^2F(\frac{\pi}{2},1,0)=0$, $\partial_pF(\frac{\pi}{2},1,0)=\partial_{pp}^2F(\frac{\pi}{2},1,0)=\partial_{qq}^2F(\frac{\pi}{2},1,0)=2$, so from previous formula the approximation is

\begin{equation}
	y(x)= 1 - \frac{1}{2}[y'(x)]^2
\end{equation}

\noindent The solution to this non linear equation, with condition $y(\frac{\pi}{2})=1$, is (computing it using {\ttfamily DSolve[\{y[x]==1-1/2*y'[x]\^{}2,y[Pi/2]==1\},y[x],x]} of {\it Mathematica})

\begin{equation}
	u(x) = 1 - \frac{\pi^2}{8} + \frac{\pi}{2}x - \frac{1}{2}x^2
\end{equation}

\noindent This expression is the second order truncated Taylor expansion of $sin(x)$ centered at $x=\frac{\pi}{2}$, so we have obtained a good approximation for the local solution found in Example 2. Note that, being in this case $\partial_xF=\partial_{xx}^2F=0$, from (\ref{taylorx}) and (\ref{taylorx2}) follows that $\partial_{xx}^2\psi(\frac{\pi}{2},0)=0$, hence second order derivative respect $x$-variable is not useful for a more accurate expression of the local solution.

\section{Relationship with series expansion of solution}

The last example seems establish a relation between symbolic expression of local solution found using implicit function technique and powers expression of solution found using series expansion of the unknown function in the original ODE. Let us consider the following symbolic expansion of second order, centered at $x_0=\frac{\pi}{2}$, for the first member of the original equation $y(x)^2+y'(x)^2-1=0$ (using the {\it Mathematica} command {\ttfamily ode = Series[y[x]\^{}2 + (D[y[x],x])\^{}2 - 1, \{x, Pi/2, 2\}]}, see \cite{diffEqMath}):

\begin{eqnarray}
&&-1+y(\frac{\pi}{2})^2+y'(\frac{\pi}{2})^2-y'(\frac{\pi}{2})\left[y(\frac{\pi}{2})+y''(\frac{\pi}{2})\right](\pi-2x)+ \nonumber \\
&&\left[y'(\frac{\pi}{2})^2+y''(\frac{\pi}{2})\left(y(\frac{\pi}{2})+y''(\frac{\pi}{2})\right)+y'(\frac{\pi}{2})y'''(\frac{\pi}{2})\right]\left(-\frac{\pi}{2}+x\right)^2 \nonumber
\end{eqnarray}

\noindent Then we solve the algebraic system formed by the expanded differential equation and the two conditions $y(\frac{\pi}{2})=1$, $y'(\frac{\pi}{2})=0$ (command {\ttfamily Solve[\{ode == 0, y[Pi/2] == 1, (D[y[x],x] /. x->Pi/2) == 0\}]}). The result consists of two possible choices: $y''(\frac{\pi}{2})=0$ and $y''(\frac{\pi}{2})=-1$. The expansion of the solution is $y(x)=y(\frac{\pi}{2})+y'(\frac{\pi}{2})(x-\frac{\pi}{2})+\frac{1}{2}y''(\frac{\pi}{2})(x-\frac{\pi}{2})^2$; for $y''(\frac{\pi}{2})=0$ we obtain the trivial solution $y=1$, for $y''(\frac{\pi}{2})=-1$ the solution is $y(x) = 1 - \frac{\pi^2}{8} + \frac{\pi}{2}x - \frac{1}{2}x^2$.\\
Therefore the series expansion method is computationally equivalent to the method of implicit function, but not logically: the first is based on a direct expansion of the unknown function in the original differential equation, while the second is based on an approximation of this ODE for finding a symbolic solution of the new equation.

\section{Finer examples}

\noindent {\it Example 3}. Let us consider the non linear ODE (see \cite{sanchez})

\begin{equation}\label{ode3}
	y'(x) = -3sin(x)\left[y(x)\right]^{\frac{4}{3}}
\end{equation}

\noindent This equation is integrable by elementary functions:

\begin{equation}\label{genint}
	y(x) = - \frac{27}{\left[c + 3cos(x)\right]^3}
\end{equation}

\noindent The associated function $F$ is $F(x,p,q)=-3sin(x)p^{\frac{4}{3}}-q$, and consider the point $T_0=(\frac{\pi}{3},1,-3\frac{\sqrt{3}}{2})$ where $F$ is null. The following identities hold:

\begin{eqnarray}
		\partial_xF=-3cos(x)p^{\frac{4}{3}}, \hspace{0.5cm} \partial_pF=-4sin(x)p^{\frac{1}{3}}, \hspace{0.5cm} \partial_qF(T_0)=-1 \\
	\partial_xF(T_0)=-\frac{3}{2}, \hspace{0.5cm} \partial_pF(T_0)=-2\sqrt{3}, \hspace{0.5cm} \partial_qF(T_0)=-1
\end{eqnarray}

\noindent We want a Taylor expansion for the implicit function $q=\phi(x,p)$ of order $(3,1)$, that is of third order respect to $x$-variable and first order respect to $p$-variable. We use expressions (\ref{taylorx}) and (\ref{taylorx2}) for $\phi$, and the following formula for the third derivative:

\begin{eqnarray}\label{taylorx3}
	\partial_{xxx}^3\phi(x_0,p_0) = -\frac{\left[3\partial_x\phi\partial_{xx}^2\phi\partial_{pp}^2F + \left[\partial_x\phi\right]^3\partial_{ppp}^3F + 3\partial_{xx}^2\phi\partial_{xp}^2F\right]}{\partial_qF} - \\ \nonumber{-\frac{\left[3\left[\partial_x\phi\right]^2\partial_{xpp}^3F + 3\partial_x\phi\partial_{xxp}^3F + \partial_{xxx}^3F\right]}{\partial_qF}}
\end{eqnarray}

\noindent Computation gives the following values for the derivatives at $T_0$:

\begin{equation}
	\partial_x\phi=-\frac{3}{2}, \hspace{0.5cm} \partial_{xx}^2\phi=3\frac{\sqrt{3}}{2}, \hspace{0.5cm} \partial_{xxx}^3\phi=3\frac{\sqrt{3}}{2}, \hspace{0.5cm} \partial_p\phi=-2\sqrt{3}
\end{equation}

\noindent so the new approximated ODE is

\begin{eqnarray}
	y'(x) = -3\frac{\sqrt{3}}{2} - \frac{3}{2}\left(x-\frac{\pi}{3}\right) + \frac{3}{4}\sqrt{3}\left(x-\frac{\pi}{3}\right)^2 \\
	\nonumber{ + \frac{\sqrt{3}}{4}\left(x-\frac{\pi}{3}\right)^3 - 2\sqrt{3}\left(y(x)-1\right)}
\end{eqnarray}

\noindent Solving this equation with conditions $y(\frac{\pi}{3})=1$, $y'(\frac{\pi}{3})=-3\frac{\sqrt{3}}{2}$ (e.g with command {\ttfamily DSolve}) one obtains a solution of the form $y(x)=e^{ax}\left(b+cx+dx^2+ex^3\right)$ which Taylor expansion of third order centered at $\frac{\pi}{3}$ is

\begin{equation}
	y(x) = 1 + \frac{\sqrt{3}}{2}\left(\pi-3x\right) + \frac{5}{12}\left(\pi-3x\right)^2 + \frac{1}{4\sqrt{3}}\left(\pi-3x\right)^3
\end{equation}

\noindent The exact solution of the original differential equation is $y(x) = - \frac{27}{\left[-\frac{9}{2} + 3cos(x)\right]^3}$ which Taylor expansion of third order centered at $\frac{\pi}{3}$ is

\begin{equation}
	y(x) = 1 + \frac{\sqrt{3}}{2}\left(\pi-3x\right) + \frac{5}{12}\left(\pi-3x\right)^2 + \frac{2}{9\sqrt{3}}\left(\pi-3x\right)^3
\end{equation}

\noindent Therefore the implicit function method gives a very good approximated symbolic solution. \\

\begin{figure}[ht]
	\begin{center}
	\includegraphics[width=8cm]{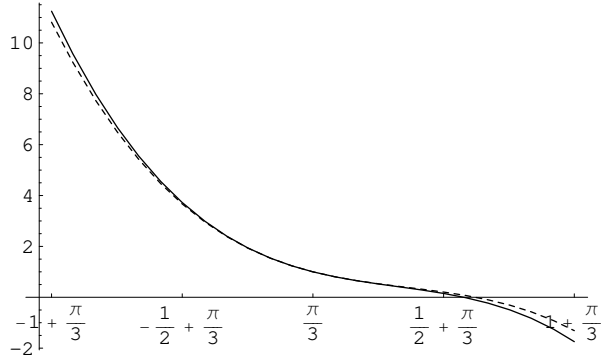}
	\caption{\small{\it Comparison between exact and approximated solution in the interval $\left[\frac{\pi}{3}-1,\frac{\pi}{3}+1\right]$; the exact solution is dashed.}} 
	\end{center}
\end{figure}

\newpage

\noindent {\it Example 4}. Consider a spherical bubble of a compressible gas immersed in an incompressible nonviscous fluid. If this fluid is flowing in a vane of a mechanical structure as rotating pump, presence of bubbles ({\it cavitation}) can give increment of noise and reduction of efficiency. Hence it is important to estimate the time $t_{c}$ necessary for the collapse of the bubble, assuming that gas pressure is neglegible compared to external fluid pressure $p_{f}$. If $y(t)$ is the radius of the bubble, it is possible to derive from Navier-Stokes equations the following ODE (see \cite{madani})

\begin{equation}\label{ode4temp}
	\frac{2}{3}y(t)y''(t) + y'(t)^2 = -\frac{2}{3}\frac{p_f}{\rho}
\end{equation}

\noindent It is convenient to express this ODE as a first order equation, which is equivalent if one doesn't consider constant solutions $y(t)=k$ (see \cite{madani} for details):

\begin{equation}\label{ode4}
	y(t)^3y'(t)^2=\frac{2}{3}\frac{p_f}{\rho}\left[R_0^3-y(t)^3\right]
\end{equation}

\noindent where $\rho$ is the density of the external fluid and $R_0$ is the radius of the bubble at time $t=0$. Without considering the trivial solution $y(t)=R_0$, which is not physically admissible, this equation is not symbolically integrable using elementary functions.\\
For simplicity, but without lack of generality, suppose $\rho=1$ and the pressures normalized with respect to the value of $p_f$, so that $p_f=1$. Let then $F(t,p,q) = \frac{2}{3}R_0^3 - \left(\frac{2}{3} + q^2\right)p^3$ the associated function $F$ defining the implicit solution. From $F(0,R_0,q_0)=0$ follows $q_0=0$, as required from the physics of the problem, therefore the point $(t_0,p_0,q_0)$ is $T_0=(0,R_0,0)$. The following formulas hold for first and second derivatives of $F$:

\begin{eqnarray}
	\nonumber{\partial_tF=0, \hspace{0.3cm} \partial_pF=-(2+3q^2)p^2, \hspace{0.3cm} \partial_qF=-2p^3q,} \\
	\nonumber{\partial_{pp}^2F=-(4+6q^2)p, \hspace{0.3cm} \partial_{qq}^2F=-2p^3, \hspace{0.3cm} \partial_{pq}^2F=-6p^2q}
\end{eqnarray}

\noindent from which we have $\partial_pF(T_0)=-2R_0^2$, $\partial_qF(T_0)=\partial_{pq}^2F(T_0)=0$, $\partial_{pp}^2F(T_0)=-4R_0$ and $\partial_{qq}^2F(T_0)=-2R_0^3$. Then, using (\ref{taylor1}), (\ref{taylorx}) and (\ref{taylorq2}), we can use the following approximation of order (1,2) (1 for $y$, 2 for $y'$), centered at $(0,R_0,0)$, for (\ref{ode4}):

\begin{equation}\label{approxode4}
	y(t) = R_0\left[1 - \frac{1}{2}y'(t)^2\right]
\end{equation}

\noindent This equation is a non linear ODE, but it is integrable by elementary functions. Its particular solution with initial condition $y(t)=R_0$ is (e.g using command {\ttfamily DSolve[\{y[t]==R\_0(1-1/2 y'[t]\^{}2),y[0]==R\_0\},y[t],t]})

\begin{equation}
	y(t) = -\frac{1}{2R_0}t^2 + R_0
\end{equation}

\noindent from which we have our estimate for the time $t_c$ of collapse:
	
\begin{equation}
	t_c = \sqrt{2}R_0
\end{equation}

\noindent From this formula one can observe the importance of the initial radius of the bubble and consequently develop suitable actions. It is interesting to compare this result with the estimate from a numerical resolution of the original ode (\ref{ode4temp}) ({\ttfamily NDSolve[\{2/3y[t]y''[t]+y'[t]\^{}2==-2/3,y[0]==R\_0\},y[t],\{t,0,0.12\}]}).

\begin{figure}[ht]\label{comparison2}
	\begin{center}
	\includegraphics[width=8cm]{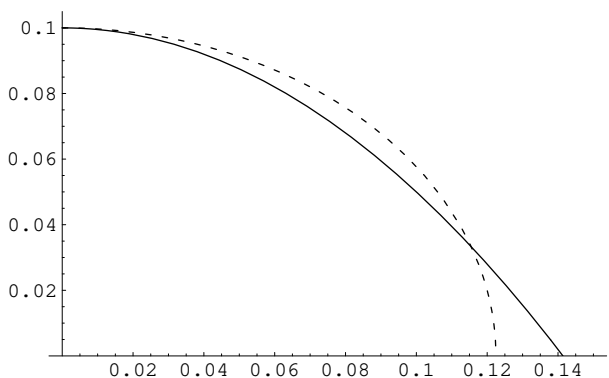}
	\caption{\small{\it Comparison between symbolic approximated solution of the reduce ODE (continuous line) and numerical solution of the original ODE (dashed line); case {\normalfont $p_{f}=1$, $\rho=1$, $R_0=0.1$}; time is in horizontal axis, radius in vertical axis.}} 
	\end{center}
\end{figure}

\noindent The figure 2 shows that symbolic solution of the reduce ODE is a good approximation of the implicit solution of the original equation.

\section{Technical notes}

We have use internal commands {\ttfamily DSolve} and {\ttfamily NDSolve} of {\it Mathematica} ver. $5.2$ for symbolic and numerical resolution of ODEs cited in the examples. Integrability of these equations using elementary functions has been tested using {\ttfamily DSolve} (see \cite{diffEqMath}).


\begin{thebibliography}{9}

\bibitem{boyce} Boyce, W.E. and DiPrima, R.C., {\it Elementary Differential Equations}, 6th edition, Wiley, (1996)

\bibitem{brota} Birkhoff, G. and Rota G., {\it Ordinary Differential Equation}, 4th edition, Wiley, New York, (1989)

\bibitem{buttazzo} Buttazzo, G. and Kawohl, B., {\it On Newton's Problem of Minimal Resistence}, in Wilson, R. and Gray, J., {\it Mathematical Conversations}, Springer-Verlag, New York, (2001)

\bibitem{courant} Courant, R. {\it Differential and Integral Calculus}, Vol. 2, Wiley, (1988)

\bibitem{diffEqMath} Coombes, K.; Hunt, B.; Lipsman, R.; Osborn, J. and Stuck, G. {\it Differential Equations with Mathematica}, 2nd edition, John Wiley \& Sons, Inc. (1998)

\bibitem{madani} Malek-Madani, R. {\it Advanced Engineering Mathematics}, Addison-Wesley, (1998)

\bibitem{mingari} Mingari, G. and Ritelli, D. {\it An Historical Outline of the Theorem of Implicit Function}, Divulgaciones Matematicas {\bf 10}(2), 171-180, (2002)

\bibitem{sanchez} Sanchez, D. {\it Ordinary Differential Equations and Stability Theory: an introduction}, Dover Publications, New York, (1979)

\bibitem{valent} Valent, T. {\it Boundary value problems of finite elasticity: local theorems on existence, uniqueness, and analytic dependence on data}, Springer, New York, (1988)

\end{thebibliography}
\end{document}